\documentclass{amsart}

\usepackage{latexsym}
\usepackage{amsfonts}
   \usepackage{amsmath}
   \usepackage{amssymb}
   \usepackage{amscd}
   \usepackage{amstext}
   \usepackage{amsthm}

\thanks{
2000 Mathematics Subject Classification 46L05(primary);
    46L08(secondary).\\
Research partially supported by
Ritsumeikan University's Fund
for sabbatical leave from the University to conduct research abroad
and JSPS Grant for Scientific Research No.14540217(c)(1).}

\begin{document}
\title{Topological stable rank of inclusions of unital C*-algebras}

\newcommand{\tsr}{\rm tsr}

\begin{abstract}
Let $1 \in A \subset B$ be an inclusion of C*-algebras
of C*-index-finite type with depth 2. We try to compute
topological stable rank of $B$ ($= \tsr(B)$)
when $A$ has topological stable rank one.
We show that $\tsr(B) \leq 2$ when $A$ is a tsr boundedly divisible
algebra, in particular,  $A$ is a C*-minimal tensor product $UHF \otimes
 D$
with $\tsr(D) = 1$. When $G$ is a finite group and $\alpha$ is an action of
$G$ on UHF, we know that a crossed product algebra
$UHF \rtimes_\alpha G$ has topological stable rank less than or equal to
two.

These results are affirmative datum
to a generalization of a  question by B. Blackadar in 1988.
\end{abstract}
\author[Hiroyuki Osaka and Tamotsu Teruya]
{Hiroyuki Osaka and Tamotsu Teruya} 
\vskip 3mm
\dedicatory{Dedicated to Jun Tomiyama on his 70th birthday}
\date{}
\address{}
\maketitle

\newtheorem{Th}{Theorem}[section]
\newtheorem{pro}[Th]{Proposition}
\newtheorem{coro}[Th]{Corollary}
\newtheorem{ho}[Th]{Lemma}
\newtheorem{rk}[Th]{Remark}
\newtheorem{df}[Th]{Definition}
\newtheorem{ex}[Th]{Example}
\newtheorem{qs}[Th]{Question}

\newcommand{\ind}{\rm Index}
\newcommand{\qe}{\square}



\section{Introduction}
The notion of topological stable rank for a C*-algebra $A$,
 denoted by $\tsr(A)$, was introduced by M. Rieffel, which
generalizes the concept of dimension of a topological space
(\cite{rf1}).
He presented the basic properties and stability theorem related to
K-Theory for C*-algebras. In \cite{rf1} he proved that
$\tsr(A \times_\alpha {\mathbb Z}) \leq \tsr(A) + 1$, and asked if
an irrational rotation algebra $A_\theta$ has topological stable rank two.
I. Putnum (\cite{pu}) gave a complete answer to this question, that is,
$\tsr(A_\theta) = 1$. Moreover, using the notion of approximate
divisibility and U. Haggerup's striking result (\cite{ha}),
B. Blackadar, A. Kumjian, and M. R\o rdam (\cite{bkr}) proved that
every nonrational
 noncommutative torus has topological stable rank one.
Naturally, we pose a question of the how to compute topological
stable rank of $A \times_\alpha G$ for any discrete group $G$.

On the contrary, one of long standing problems is whether a fixed point
algebra of a UHF C*-algebra by an action of a finite group $G$ is an
AF C*-algebra.
O. Bratteli (\cite{obr})  proved that any fixed point algebra of an
UHF-algebra by a product type action of a finite abelian group
is an AF C*-algebra. In 1988, B. Blackadar (\cite{bl3}) constructed a
symmetry
on the CAR algebra whose fixed point algebra is not an AF C*-algebra.
Note that A. Kumjian (\cite{km}) constructed a symmetry
on a simple AF C*-algebra
whose fixed point algebra is not an AF C*-algebra.
Later, D. Evans and A. Kishimoto proved that for any compact group
$G \not= \{e\}$ and $p \geq 2$, there exists an action of $G$ on
$M_{p^\infty}$ whose fixed point algebra is not an AF C*-algebra.
All these constructions embodied expressing the AF C*-algebra $A$ as an
inductive limit $A_1 \rightarrow A_2 \rightarrow \cdots \rightarrow$
$A = \lim A_n$, where
each C*-algebra $A_n$ is not an AF C*-algebra.
This is related to the classification theory of simple unital AH-algebras
(\cite{dg},\cite{ell0},\cite{eg}).
Indeed, applying the classification theory
G. Elliott constructed a symmetry $\alpha$ on an UHF algebra, and
proved that $UHF \rtimes_\alpha({\mathbb Z}/2{\mathbb Z})$
is not AF C*-algebra,
but AI-algebra, that is, the inductive limit of direct sums of
$C([0,1]) \otimes M_n({\mathbb C})$. Note that this crossed product algebra
has topological stable rank one and real rank one. B. Blackadar proposed
the question in the same article whether $\tsr(A \times_\alpha G) = 1$
for any unital AF C*-algebra $A$,
a finite group $G$, and an action $\alpha$ of $G$ on $A$.

In this paper we try
to solve B. Blackadar's question from more general situation
using C*-index theory defined by Y. Watatani (\cite{wata}).
In the case that an inclusion $1 \in A \subset B$ is of
index-finite type with depth 2 if $A$ is tsr boundedly divisible algebra
(see \cite[Definition 4.1]{rf3} and section 5) with $\tsr(A) = 1$
we show that $\tsr(B) \leq 2$ (Theorem 5.1). Hence if $A$ is
a UHF C*-algebra, we conclude that $\tsr(B) \leq 2$
under the above condition.
Therefore we get an affirmative data to B. Blackadar's question.
Namely, for any UHF C*-algebra $A$, a finite group $G$, and an
action $\alpha$ of $G$ on $A$, we conclude that 
a crossed product algebra $A \times_\alpha G$
has $\tsr(A \times_\alpha G) \leq 2$.
We can not still get the complete answer, but it seems to guarantee that
the question  would be solved affirmatively.

This paper is organized as follows.
In section 2 we state a number of preliminary results about
topological stable rank.
In section 3 we explain a brief survey of C*-index theory.
In section 4 we study the quasi-basis for the induced  conditional
expectation
of the derived  inclusion $1 \in pAp \subset pBp$
from the inclusion $1 \in
A \subset B$ and a non-zero projection $p \in A$.
We give a new estimate of topological stable rank for the
inclusion of index-finite type with depth 2 in section 5 and give the
main theorem, that is, that
topological stable rank of a crossed product algebra
of a UHF algebra by any finite group and  any action
has less than or equal to 2.

\vskip 5mm

\leftline{\bf Acknowledgement}

Some results of this note were conducted while the first author was
in sabbatical leave visiting the University of Oregon.
He would like to thank the member of
the mathematics department there for their warm hospitality. In particular,
he would like to thank Huaxin Lin and N. Christopher Phillips for many
stimulating conservations about ranks.
He also would like to thank Ken Goodearl for his fruitful discussion about
cancellation property.

The second author would like to thank the Principal Masajiro Nashiro of
Okinawa Shogaku  high school
for his moral and material support.

\vskip 3mm

\section{Topological stable rank}

In this section we present a definition of topological stable rank
and some basic estimate
formulas for it.

\begin{df}
Let $A$ be a unital C*-algebra and $Lg_n(A)$ be
the set of elements $(b_i)$ of $A^n$ such that
$$
Ab_1 + Ab_2 + \cdots + Ab_n = A.
$$
 Then topological stable rank of $A$,
$\tsr(A)$, to
be the least integer $n$ such that the set $Lg_n(A)$ is dense in $A^n$.
Topological stable rank of a non-unital C*-algebra is defined by
topological  stable rank of its unitaization algebra $\tilde{A}$
\end{df}

Note that $\tsr(A) = 1$ is equivalent to have the dense set of
invertible elements in $\tilde{A}$.  Here are some formulas for computing
stable rank of C*-algebras.

\begin{ho}
Let $A$ be a unital C*-algebra, and let $a_1, a_2, \dots, a_n \in A$.
The followings are equivalent:

\vskip 5mm

\begin{enumerate}
\item[(1)]
There are $c_1, c_2, \dots, c_n \in A$ such that
$c_1a_1 + c_2a_2 + \dots + c_na_n$ is invertible.
\item[(2)]
There are $c_1, c_2, \dots, c_n \in A$ such that
$c_1a_1 + c_2a_2 + \cdots + c_na_n = 1$.
\item[(3)]
$a_1^*a_1 + a_2^*a_2 + \cdots + a_n^*a_n$ is invertible.
\end{enumerate}
\end{ho}

\vskip 3mm

{\it Proof}. Standard. \hfill$\Box$

\vskip 5mm

\begin{Th}[\cite{rf1}]
\begin{enumerate}
\item[(1)]
Let $J$ be a closed two-sided ideal of C*-algebra $A$.
Then
$$
\tsr(A) \leq \max\{\tsr(J), \tsr(A/J) + 1\}.
$$
\item[(2)]
Let $n$ be a positive integer. Then
$$
\tsr(M_n(A)) = \left\{\frac{\tsr(A) -1}{n}\right\} + 1,
$$
where $\{t\}$ denotes
the least integer which is greater than or equal to $t$.
\item[(3)]
Let $p$ be a non-zero projection in $A$.
Then $\tsr(A) = 1$ if and only if $\tsr(pAp) = \tsr((1- p)A(1- p)) = 1$.
\item[(4)]
Let $\alpha$ be an action of $A$. Then
$$
\tsr(A \rtimes_\alpha {\mathbb Z}) \leq \tsr(A) + 1.
$$
\item[(5)]
Let ${\mathbb K}$ be a C*-algebra of compact operators on an
 infinite dimensional separable Hilbert space. Then
$$
\tsr(A) = 1 \ \hbox{if and only if} \ \tsr(A \otimes {\mathbb K}) = 1.
$$
\end{enumerate}
\end{Th}

\vskip 5mm

{}From the point of C*-module we
 have the following
formula for topological stable rank.

\begin{Th}
Let $1 \in A \subset B$ be an inclusion of C*-algebras.
Suppose that there are elements $\{v_i\}_{i=1}^n \in B$
such that
$$
B = Av_1 + Av_2 + \cdots + Av_n.
$$
Then
$$
\tsr(B) \leq n \times \tsr(A).
$$
\end{Th}

\vskip 3mm

{\it Proof}.
The proof is the same as in \cite[Theorem~ 2.1]{jo}.
We will put a sketch of its proof for a self-contained.

We first assume that $\tsr(A)=1$. Let $\{a_{i1}v_1 + \cdots +
a_{in}v_n\mid i=1,\dots,n\}$ be $n$ elements in $B$.
Then,
$$
\begin{array}{ll}
\left(\begin{array}{c}
       a_{11}v_1+a_{12}v_2+\cdots+a_{1n}v_n\\
       a_{21}v_1+a_{22}v_2+\cdots+a_{2n}v_n\\
       \vdots\\
       a_{n1}v_1+a_{n2}v_2+\cdots+a_{nn}v_n
       \end{array}
       \right)
&= \left(\begin{array}{cccc}
      a_{11}&a_{12}&\cdots&a_{1n}\\
      a_{21}&a_{22}&\cdots&a_{2n}\\
      \vdots&\vdots&\ddots&\vdots\\
      a_{n1}&a_{n2}&\cdots&a_{nn}

      \end{array}
      \right)
\left(\begin{array}{c}
     v_1\\
     v_2\\
     \vdots\\
     v_n
     \end{array}
     \right)\\
&= av
\end{array}
$$
Since $\tsr(M_n(A))=1$ by Theorem~2.3 (2),
we can approximate the matrix $(a_{ij})$
close enough by an invertible  matrix $x=(x_{ij})$  in
$M_n(A)$.  Then the element $y = xv$  is close to $av$.
Since $(v_1, v_2, \dots, v_n) \in Lg_n(B)$,
$(x_{11}v_1+\cdots +x_{1n}v_n,\dots ,x_{n1}v_1+\cdots+x_{nn}v_n)$
belongs to $Lg_n(B)$ by Lemma~ 2.2, and  is close to
$(a_{11}v_1+\cdots +a_{1n}v_n,\dots ,a_{n1}v_1+\cdots+a_{nn}v_n)$, which
completes the proof in the case that $\tsr(A)=1$.

Similarly, one can prove  the theorem when
$\tsr(A)>1$ using Theorem~ 2.3 (2).
\hfill$\qe$

\vskip 5mm

\begin{coro}
Let $1 \in A \subset B$ be an inclusion of C*-algebras, and
let $E$ be a faithful conditional expectation from $B$ to $A$
of index finite type.
That is,
there is a quasi-basis $\{(v_i, v_i^*)
\}_{i=1}^n$  in $B \times B$
such that
any element $x \in B$ can be written as
$$
x = \sum_{i=1}^nE(xv_i)v_i^* = \sum_{i=1}^nv_iE(v_i^*x).
$$
Then
$$
\tsr(B) \leq n \times \tsr(A).
$$
\end{coro}

{\it Proof}.
Since
$$
B = Av_1^* + Av_2^* + \cdots + Av_n^*,
$$
we are done by Theorem~ 2.4.
\hfill$\Box$

\vskip 5mm

\begin{coro}
Let $A$ be a unital C*-algebra and let $G$ be a finite group.
Then
$$
\tsr(A \rtimes_\alpha G) \leq |G|\tsr(A),
$$
where $|G|$ is the cardinal number of $G$.
\end{coro}

{\it Proof}.
Let $\alpha \colon G \rightarrow Aut(A)$ be a representation and we
assume that the crossed product $A \rtimes_\alpha G$ acts on some Hilbert
space. Let $\{u_g\}_{g \in G}$ be implemented unitaries of $\alpha_g$
such that $\alpha_g(a) = u_gau_g^*$, for all $a \in A$.
Then any element $x$ in
$A \rtimes_\alpha G$ can be written as
$x = \sum_{g \in G}a_gu_g$.
Let $E \colon A \rtimes_\alpha G \rightarrow A$ be the
canonical conditional expectation by $E(x) = a_e$.
Then,
$$
x = \sum_{g \in G}E(xu_g^*)u_g, \quad \forall x \in A \rtimes_\alpha G.
$$
Therefore it follows from Corollary~ 2.5 that
$$
\tsr(A \rtimes_\alpha G) \leq |G|\tsr(A).
$$
\hfill$\qe$

\vskip 5mm

\begin{rk}
The estimate of topological stable rank of crossed product in
{\em Corollary~ 2.6} is not best one.
Indeed, in the case of $G = {\mathbb Z}/n{\mathbb Z}$
we have
$$
\tsr(A \rtimes_\alpha G) \leq \tsr(A) + 1.
$$
using {\em Theorem~ 2.3 (4)}.
But, in section 5 we show more better estimate
in the case that $A$ is tsr boundedly divisible
with $\tsr(A) = 1$, $G$ any finite group, and $\alpha$
an action of $G$ on $A$  as follows:
$$
\tsr(A \rtimes_\alpha G) \leq 2.
$$
\end{rk}

\vskip 5mm

\section{C*-Index Theory}
In this section we summarize the C*-index theory of Y. Watatani
(\cite{wata}).

Let $1 \in A \subset B$ be an inclusion  of C*-algebras, and let $E\colon B
\rightarrow A$ be a faithful conditional expectation from $B$ to $A$.

A finite family
$\{(u_1, v_1), \dots, (u_n, v_n)\}$ in $B \times B$ is called {\em a

quasi-basis} for $E$ if $$
\sum_{i=1}^n u_iE(v_ib) = \sum_{i=1}^nE(bu_i)v_i = b \enskip \hbox{for}
\enskip b \in B.
$$
We say that a conditional expectation $E$ is  of {\em index-finite type}
if there
exists a quasi-basis for $E$. In this case the index of $E$ is defined by
$$
{\rm Index}(E) = \sum_{i=1}^nu_iv_i.
$$
(We say also  that the inclusion $1 \in A \subset B$ is
of {\em index-finite type}.)

Note that ${\rm Index}(E)$ does not depend on the choice of a quasi-basis
(\cite[Example~ 3.14]{iz3})
and every conditional expectation $E$ of index-finite type on a C*-algebra
has a quasi-basis of the form $\{(u_1, u_1^*), \dots, (u_n,u_n^*)\}$
(\cite[Lemma~ 2.1.6]{wata}).
Moreover ${\rm Index}(E)$ is always contained in the center of $B$, so
that it is a scalar whenever $B$ has the trivial center, in particular when
$B$ is simple (\cite[Proposition~ 2.3.4]{wata}).

Let $E\colon B\to A$ be a faithful conditional expectation.
Then $B_{A}(=B)$ is
 a pre-Hilbert module over $A$ with an $A$-valued inner
product $$\langle x,y\rangle =E(x^{*}y), \ \ x, y \in B_{A}.$$
Let $\mathcal E$ be the completion of $B_{A}$ with respect to the norm on
$B_{A}$ defined by
$$\| x\|_{B_{A}}=\|E(x^{*}x)\|_{A}^{1/2}, \ \ x \in B_{A}.$$
Then $\mathcal E$ is a Hilbert $C^{*}$-module over $A$.
Since $E$ is faithful, the canonical map $B\to \mathcal E$ is injective.
Let $L_{A}(\mathcal E)$ be the set of all (right) $A$-module homomorphisms
$T\colon \mathcal E\to \mathcal E$ with an adjoint $A$-module homomorphism
$T^{*}\colon \mathcal E\to \mathcal E$ such that $$\langle T\xi,\zeta
\rangle =
\langle \xi,T^{*}\zeta \rangle \ \ \ \xi, \zeta \in \mathcal E.$$
Then $L_{A}(\mathcal E)$ is a $C^{*}$-algebra with the operator norm
$\|T\|=\sup\{\|T\xi \|:\|\xi \|=1\}.$ There is an injective
$*$-homomorphism $\lambda \colon B\to L_{A}(\mathcal E)$ defined by
$$\lambda(b)x=bx$$
for $x\in B_{A}$ and  $b\in B$, so that $B$ can
be viewed as a
$C^{*}$-subalgebra of $L_{A}(\mathcal E)$.
Note that the map $e_{A}\colon B_{A}\to B_{A}$
defined by $$e_{A}x=E(x),\ \ x\in
B_{A}$$ is
bounded and thus it can be extended to a bounded linear operator, denoted
by $e_{A}$ again, on $\mathcal E$.
Then $e_{A}\in L_{A}({\mathcal E})$ and $e_{A}=e_{A}^{2}=e_{A}^{*}$; that
is, $e_{A}$ is a projection in $L_{A}(\mathcal E)$.
A projection $e_A$ is called the {\em Jones projection} of $E$.

The {\sl (reduced) $C^{*}$-basic construction} is a $C^{*}$-subalgebra of
$L_{A}(\mathcal E)$ defined to be
$$
C^{*}(B, e_{A}) = \overline{ span \{\lambda (x)e_{A} \lambda (y) \in
L_{A}({\mathcal E}): x, \ y \in B \ \} }^{\|\cdot \|} $$
(\cite[Definition~ 2.1.2]{wata}).

Then we have

\begin{ho}{\rm (\cite[Lemma 2.1.4]{wata}) } 
\begin{enumerate}
\item[$(1)$]
$e_{A}C^{*}(B, e_{A})e_{A}=
\lambda(A)e_{A}.$
\item[$(2)$]
 $\psi :A\to e_{A}C^{*}(B, e_{A})e_{A}$, $\psi (a)=\lambda (a)e_{A}$, is
a $*$-isomorphism (onto).
\end{enumerate}
\end{ho}

\begin{ho}{\rm (\cite[Lemma~ 2.1.5]{wata})} The following are equivalent:

\begin{enumerate}
\item[$(1)$] 
$E:B\to A$ is of index-finite type.
\item[$(2)$]
 $C^{*}(B, e_{A})$ has an identity and there exists a number $c$ with 
 $0 < c < 1$ such that 
$$
E(x^{*}x)\geq c(x^{*}x) \ \ \ x\in B.
$$
\end{enumerate}
\end{ho}

The above inequality was shown first in \cite{pp} by Pimsner and Popa for
the conditional expectation $E_{N}\colon M\to N$
from a type II$_{1}$ factor $M$
onto its subfactor $N$ ($c$ can be taken as the inverse of the Jones index
$[M : N]$).

A conditional expectation $E_{B}\colon C^{*}(B, e_{A})\to B$ defined by
$$
E_{B}(\lambda (x)e_{A}\lambda (y))=({\rm Index}(E))^{-1}xy \ \hbox{for}
\ x \ \hbox{and}\ y \in B
$$
is called {\sl the dual conditional expectation} of $E:B \to A$. If $E$ is
of index-finite type, so is $E_B$ with a quasi-basis $\{(w_i, w_i^*)\}$,
where $w_i = u_ie_A\ind(E)^{\frac{1}{2}}$, and $\{(u_i, u_i^*)\}$ is a
quasi-basis for $E$ (\cite[Proposition~ 2.3.4]{wata}).

Even if ${\rm Index}(E)$ is scalar, we do not know the relation between the
number of pairs in a quasi-basis and ${\rm Index}(E)$
(\cite[Example~ 3.14]{iz3}\cite[Lemma~ 3.4]{o}).
Izumi
, however, showed recently that if we extend a conditional expectation
$E$ from $\sigma$-unital C*-algebra $D$ onto stable simple C*-algebra $C$
with $\overline{DC} = D$ to
a faithful conditional expectation $\tilde{E}$ from
the multiplier algebra $M(D)$ onto $M(C)$, then it has only
one pair as a quasi-basis(\cite[Proposition~ 3.6]{iz3}).

The inclusion $1 \in A \subset B$ of unital C*-algebras of
index-finite type is said to
have {\it finite depth k}
if the derived tower obtained by iterating the basic construction
$$
A' \cap A \subset A' \cap B \subset
A' \cap B_2 \subset A' \cap B_3 \subset
\cdots
$$
satisfies
$(A' \cap B_k)e_k(A' \cap B_k) = A' \cap B_{k+1}$,
where $\{e_k\}_{k \geq 1}$ are projections
derived obtained by iterating the basic construction such that
$B_{k+1} = C^*(B_{k}, e_k)$ \ ($k \geq 1$) \quad ($B_1 = B, e_1 = e_A$).
Let $E_k : B_{k+1} \rightarrow B_k$ be a faithful conditional expectation
correspondent to $e_
k$ for $k \geq 1$.
Moreover we have

\vskip 3mm
\begin{ho}{\rm (\cite[Lemma 2.3.5]{wata}) }
Suppose $\ind(E)$ is in $A$. Then
$$
\left\{
\begin{array}{l}
e_{k+1}e_ke_{k+1} = (\ind(E))^{-1}e_{k+1}\\
e_ke_{k+1}e_k = (\ind(E))^{-1}e_k
\end{array}
\right.
$$
for $1 \leq k$.
\end{ho}

\vskip 5mm

When $G$ is a finite group and $\alpha$ an
action of $G$ on $A$,
it is well known that an inclusion
$1 \in A \subset  A \rtimes_\alpha G$
is of depth 2.
We will give its proof for a self-contained.

\begin{ho}
Let $A$ be a unital C*-algebra, $G$ a finite group, and
$\alpha$ an action of $G$ on $A$.
Then an inclusion $1 \in A \subset A \rtimes_\alpha G$
is of depth 2.
\end{ho}

{\it Proof}.
Assume that $A \rtimes_\alpha G$ acts on some Hilbert space $H$, and
there are implemented unitaries $\{u_g\}_{g \in G}$ such that
$\alpha_g(a) = u_gau_g^*$ and any element $x \in A \rtimes_\alpha G$
can be written by
$x = \sum_{g\in G}a_gu_g$
for some $a_g \in A$.
Let $E\colon A \rtimes_\alpha G \rightarrow A$ be the canonical
conditional expectation by $E(\sum_{g\in G}a_gu_g) = a_e$.
Then a set $\{(u_g^*,u_g)\}_{g\in G}$ is a quasi-basis for $E$.
Note that $\ind(E) = |G|$, where $|G|$ is the cardinal number of $G$.

Let
$$
1 \in A \subset A \rtimes_\alpha G \subset B_2 \subset B_3 \subset \cdots
$$
be a sequence of basic constructions.

Claim 1: $u_he_Au_h^*$ is a projection from 
${\mathcal E}_{A \rtimes_\alpha G}$
to $Au_h$ for each $h \in G$.

{\it Proof}.
It is obvious that $u_he_Au_h^*$ is projection. Since
$$
\begin{array}{ll}
u_he_Au_h^*(\sum_ga_gu_g) &= u_he_A\sum u_h^*a_gu_g\\
&= u_hE(\alpha_{h^{-1}}(a_g)u_h^*u_g)\\
&= u_h\alpha_{h^{-1}}(a_h) = a_hu_h,
\end{array}
$$
we have the claim 1.
\hfill$\Box$

\vskip 3mm

Claim 2: For any $h \in G$ $u_he_Au_h^* \in A' \cap B_2$.

{\it Proof}.
Since for any $a \in A$
$$
\begin{array}{ll}
u_he_Au_h^*a &= u_he_Au_h
^*au_hu_h^*\\
&= u_he_A\alpha_{h^{-1}}(a)u_h^*\\
&= u_h\alpha_{h^{-1}}(a)e_Au_h^*\\
&= au_he_Au_h^*,
\end{array}
$$
we have the claim 2.
\hfill$\Box$

\vskip 3mm

Claim 3: $\{u_he_Au^*_h\}_{h \in G}$ are orthogonal projections in $B_2$.

{\it Proof}.
Trivial.

\vskip 3mm

Claim 4: $(A' \cap B_2)e_2(A' \cap B_2) =  A' \cap B_3$,
where $e_2$ is a projection correspondent to the dual conditional
expectation from $B_2$ to $B_1$.

{\it Proof}.
Note that $(A' \cap B_2)e_2(A' \cap B_2)$ is a closed two-sided ideal
of $A' \cap B_3$ (\cite[Theorem~ 4.6.3(ii)]{ghj}).
So we have only to show that this ideal contains an identity.

Since $u_ge_Au_g^* \in A' \cap B_2$, we have
$$
\begin{array}{ll}
(u_ge_Au_g^*)e_2(u_ge_Au_g^*) &= u_ge_Ae_2e_Au_g^*\\
&= \frac{1}{|G|}u_ge_Au_g^*
\in (A' \cap B_2)e_2(A' \cap B_2).
\end{array}
$$

The last equality comes from Lemma 3.2.
Hence
$$
1 = \sum_{g \in G}u_ge_Au_g^* \in (A' \cap B_2)e_2(A' \cap B_2).
$$
This means that the inclusion $A \subset A \rtimes_\alpha G$ is of depth 2.
\hfill$\Box$

\vskip 5mm

\begin{rk}
When $A$ is simple and
$\alpha$ is outer, the crossed
product $A \rtimes_\alpha G$ is simple
$($\cite{ki}$)$ and  we easily have

\begin{enumerate}
\item[(1)]
$A' \cap B_2$ is isomorphic to $\sum_{g\in G}{\mathbb C}e_g$, where 
$e_g = u_ge_Au_g^*$
\item[(2)]
$A' \cap B_3$ is isomorphic to $M_{|G|}({\mathbb C})$.
\end{enumerate}
\end{rk}

\vskip 5mm


\section{Quasi-basis for a conditional expectation}

Let $1 \in A \subset B$ be an inclusion  of unital C*-algebras and
$E : B \rightarrow A$ be a faithful conditional expectation
of index-finite type. Let $\{(v_i,v_i^*)\}_{i=1}^n$
be a quasi-basis for $E$.

The following is kindly informed by Y. Watatani (\cite{wata2})
to the first author.

\vskip 3mm

\begin{pro}

Under the above situation
if a projection $p$ in $A$ has elements $\{y_j\}_{j=1}^m$ in $A$
such that $\sum_{j=1}^my_jpy_j^* = 1$, then a set
$\{(pv_iy_jp,py_j^*v_i^*p)\}_{1\leq i\leq n, 1\leq j \leq m}$
is a quasi-basis for
a conditional expectation $
F_p = E|pBp$ from $pBp$ onto
$pAp$.
Moreover $\ind(E) = \ind(F_p)$.
\end{pro}

{\it Proof}. It follows from the direct computation.
Indeed
for any $b \in B$ we have
$$
\begin{array}{ll}
\sum_{i,j}pv_iy_jpF_p(py_j^*v_i^*ppbp)
&= \sum_{i,j}pv_iy_jpE(y_j^*v_i^*ppb)p\\
&= \sum_{i=1}^npv_i(\sum_{j=1}^my_jpy_j^*)E(v_i^*ppb)p\\
&= \sum_{i=1}^npv_iE(v_i^*pb)p = p(pb)p = pbp.
\end{array}
$$
Similarly,
$$
\sum_{i,j}F_p(pbppv_iy_jp)py_j^*x_i^*p = pbp.
$$
So $\{(pv_iy_jp,py_j^*x_i^*p)\}_{1 \leq i \leq n, 1 \leq j \leq m}$
is a quasi-basis for $F_p$ and
$$
\begin{array}{ll}
{\rm Index}(F_p) = \sum_{i,j}pv_iy_jpy_j^*x_i^*p
&= \sum_{i=1}^npv_i(\sum_{j=1}^my_jpy_j^*)v_i^*p\\
&= \sum_{i=1}^npv_iv_i^*p = ({\rm Index}(E))p.
\end{array}
$$
\hfill$\qe$

\vskip 5mm

\begin{coro}
Let $1 \in A \subset  B$ be an inclusion of unital C*-algebras and
$E\colon B \rightarrow A$ be a faithful conditional expectation
of index-finite type. Suppose that
$A$ is simple.
Then for any non-zero projection $p$ in $A$
the conditional expectation $F_p = E|pBp \colon pBp
\rightarrow pAp$ is of index-finite type.
\end{coro}

{\it Proof}.
Since $A$ is simple, there are finite elements $a_i$ such that
$\sum_ia_ipa_i^* = 1$. So the statement comes from the previous
proposition.
\hfill$\qe$

\vskip 3mm

The following result shows that the
Jones projection of $F_p$ in the previous result is
$e_Ap$.

\begin{pro}
Let $1 \in A \subset B$ be an inclusion of C*-algebras with finite index,
and let $e_A$ be the Jones projection correspondent to
a faithful conditional expectation $E\colon B \rightarrow A$.
Suppose that $A$ is simple. Then for any projection $p \in A$,
$e_Ap$ is the Jones projection for the conditional expectation
$F_p = E|pBp\colon pBp \rightarrow pAp$ and
$pC^*(B,e_A)p$ is the basic construction
for $F_p$.
\end{pro}

{\it Proof}.
For any $x \in pBp$
$$
e_Apxe_Ap = E(px)e_A
p = E(pxp)e_Ap = F_p(x)e_Ap.
$$
Since $A$ is simple, the map
$$
pAp \ni x \rightarrow xe_Ap (= xe_A) \in L( \mathcal E)
$$
is injective, where $\mathcal E$ is a Hilbert $A$-module obtained by
the basic construction. Then by \cite[Proposition~ 2.2.11]{wata}
$e_Ap$ is the Jones projection and $C^*(pBp,e_Ap)$ is the basic construction
for $F$. It is obvious that $C^*(pBp,e_Ap) = pC^*(B,e_A)p$.
\hfill$\qe$

\vskip 5mm

The following result means that the number of elements in
quasi-basis is stable under the particular situation.

\begin{pro}
Let $1 \in A \subset B$ be an inclusion of C*-algebras with
finite index.
Suppose that there is a C*-subalgebra $D$ of $A$ such that
$e_A$ is full in $D' \cap C^*(B,e_A)$.
Then there are finitely elements
$\{v_i\}_{i=1}^n$ in $D' \cap B$ such that for any non-zero
projection $p \in D$ the sets
$\{(v_i, v_i^*)\}_{i=1}^n$ and $\{(pv_i, pv_i^*)\}_{i=1}^n$ are
quasi-basis for $E$ and $F_p$, respectively.
\end{pro}

{\it Proof}.

Since $(D' \cap C^*(B,e_A))e_A(D' \cap C^*(B,e_A)) =
D' \cap C^*(B,e_A)$, there are
finitely elements $\{x_i\}_{i=1}^n$ in  $D' \cap C^*(B,e_A)$
such that
$$
\sum_{i=1}^nx_ie_Ax_i^* = 1.
$$
Using the standard argument(\cite{pp}) we can find $v_i \in B$
such that $v_ie_A = x_ie_A$. Since $E$ is faithful,
$v_i \in D' \cap B$.

Since for any $b \in B$

$$
\begin{array}{ll}
b = 1\cdot b &= \sum_{i=1}^nv_ie_Av_i^*b\\
&= \sum_{i=1}^nv_iE(v_i^*b) = \sum_{i=1}^nE(bv_i)v_i^*,
\end{array}
$$
it follows that
$\{(v_i, v_i^*)\}_{i=1}^n$ is a quasi-basis for $E$.

{}From the simple calculus we know that $\{(pv_i, pv_i^*)\}_{i=1}^n$
is a quasi-basis for $F_p$.
\hfill$\qe$

\vskip 5mm

\vskip 5mm

\section{Topological stable rank of inclusions of
C*-algebras}

In this section we prove the following main result:

\begin{Th}
Let $1 \in A \subset B$ be an inclusion of unital
C*-algebras and $E\colon B \rightarrow A$ be
a faithful conditional expectation
of index-finite type.
Suppose that the inclusion $1 \in A \subset B$
has depth 2 and
$A$ is tsr boundedly divisible with $\tsr(A) = 1$.
Then $B$ is tsr boundedly divisible.
Moreover we have $\tsr(B) \leq 2$.
\end{Th}

Recall that a C*-algebra $A$ is {\it tsr boundedly divisible}
(\cite[Definition 4.1]{rf3})
if there is a constant $K$ ($> 0$) such that for every positive integer $m$
there is an integer $n \geq m$ such that $A$ can be expressed
as $M_n(B)$ for a C*-algebra $B$ with $\tsr(B) \leq K$.
For any unital C*-algebra $A$ with $\tsr(A) = 1$
a C*-minimal tensor product algebra $A \otimes UHF$
is a typical tsr boundedly divisible algebra.

\vskip 3mm

Before giving the proof, we state a useful result by B. Blackadar.

\vskip 3mm

\begin{ho}$($\cite[Lemma A6]{bl0}$)$
Let $A$ be a unital C*-algebra and $p$ be a full projection in $A$ such
that $\sum_{i=1}^nu_ipv_i = 1$ for some elements $u_i, v_i$ in $A$.
Then $\tsr(A) \leq \tsr(pAp) + n - 1$.
\end{ho}

\vskip 3mm

\begin{rk}
Very recently, B. Blackadar sharpened the estimate in the previous lemma
$($\cite{bl4}$)$.
That is, for any unital C*-algebra $A$ and non-zero projection $p \in A$
$\tsr(A) \leq \tsr(pAp)$.
\end{rk}

\vskip 5mm

The following estimate is the converse of Corollary~ 2.5.

\begin{pro}
Let $1 \in A \subset B$ be an inclusion of unital C*-algebra
of index-finite type.
Let $\{(v_i,v_i^*)\}_{i=1}^m$ be a quasi-basis for
a faithful conditional expectation $E$ from $B$ onto $A$.
Then we have
$$
\tsr(A) \leq m^2(\tsr(B) + 1) - 2m + 1.
$$
\end{pro}

{\it Proof}.
Let $B_2$ be the basic construction derived from this inclusion
and $E_2$ be the dual conditional expectation from $B_2$ onto $B$.
Note that a set \newline
$\{(v_ie_A(\ind(E))^{\frac{1}{2}},(\ind(E))^{\frac{1}{2}}e_Av_i^*)\}_{i=1}^m
$
is the  quasi-basis for $E_2$. Hence from Corollary~ 2.5 we have
$$
\tsr(B_2) \leq m \times \tsr(B).
$$

By \cite[Lemma 3.3.4]{wata}
there is an isomorphism $\varphi : B_2 \rightarrow qM_m(A)q$

such that
$$
\phi(xe_Ay) = [E(v_i^*x)E(yv_j)]_{i,j=1}^m,
$$
where $q = [E(v_i^*v_j)]$.
Note that $q$ is a projection.

{\bf Claim}: There are $X_i, Y_i \in M_m(A)$ such that
$\sum_{i=1}^mX_iqY_i = 1.$

{\it Proof of the Claim}:
Let
$$
\begin{array}{ll}
(X_i)_{h,k} = \left\{
      \begin{array}{ll}
      E(v_k)&\ \quad\hbox{if}\ h = i\\
      0&\ \hbox{others}
      \end{array}
      \right.
      ,
&(Y_i)_{h,k} = \left\{
\begin{array}{ll}
E(v_h^*)&\ \hbox{if}\ k = i\\
0&\ \hbox{others}
\end{array}
\right.
\end{array}
$$
for each $1 \leq h, k \leq m$.

Then from a simple calculation we have
$$
\sum_{i=1}^mX_iqY_i = 1.
$$

So from Lemma 5.2 we have
$$
\tsr(M_m(A)) \leq \tsr(qM_m(A)q) + m - 1.
$$
Since $\varphi(B_2)$ is isomorphic to $qM_m(A)q$ and
$\tsr(B_2) \leq m\times\tsr(B)$, we have
$$
\tsr(M_m(A)) \leq m\times\tsr(B) + m - 1 = m(\tsr(B) + 1) - 1.
$$

Since from Theorem~ 2.3 (2) we have
$$
\frac{\tsr(A) - 1 }{m} + 1 \leq \tsr(M_m(A)) \leq m(\tsr(B) + 1) - 1,
$$
it follows that
$$
\tsr(A) \leq m^2(\tsr(B) + 1) - 2m + 1.
$$
\hfill$\qe$

\vskip 3mm

{\bf Proof of Theorem 5.1}

Let

$$
1 \in A \subset B \subset B_2 \subset B_3 \subset \cdots
$$
be the derived tower of iterating the basic construction and
$\{e_k\}_{k\geq1}$ be canonical projections such that
$B_{k+1} = C^*(B_k,e_k)$, where $e_1 = e_A$.
Since $1 \in A \subset B$ is of depth 2,
we have
$$
(A' \cap B_2)e_2(A' \cap B_2) = A' \cap B_3.
$$
Hence there exit some $n \in {\mathbb N}$
a quasi-basis $\{(u_i,u_i^*)\}_{i=1}^n$
for the conditional expectation $E_2$ from $B_2$ onto $B$
so that $u_i \in A' \cap B_2$ for $1 \leq i \leq n$
(see the proof of Proposition~ 4.4).

Since $B_2$ is stably isomorphic to $A$, we know that $\tsr(B_2) = 1$
by Theorem~ 2.3 (5). Take non-zero projection $p$ in $A$.
Since $u_i \in A' \cap B_2$ for $1 \leq i \leq n$,
a set $\{(pu_i,u_i^*p)\}_{i=1}^n$ is a quasi-basis for
the conditional expectation $F_p = E_2|pB_2p$ from $pB_2p$ onto $pBp$.
Hence from Proposition~ 5.4 we have
$$
\begin{array}{ll}
\tsr(pBp) &\leq n^2(\tsr(pB_2p) +
 1) - 2n + 1\\
&= 2n^2 - 2n + 1.
\end{array}
$$
The last equality comes from that $\tsr(B_2) = 1$ and Theorem~ 2.3 (3).

Since $A$ is tsr boundedly divisible, for any $l \in {\mathbb N}$
there are $k \in {\mathbb N}$ with $k \geq l$ and a C*-algebra $D$
such that $A \cong M_k(D)$. Hence there is a matrix system
$\{e_{i,j}\}_{i,j=1}^k$ in $A$ such that
$B \cong M_k(e_{1,1}Be_{1,1})$.
Then from the above estimate we have
$$
\tsr(e_{1,1}Be_{1,1}) \leq n^2 + (n - 1)^2.
$$
Therefore $B$ is tsr boundedly divisible.
{}From Theorem~ 2.3 (2) we can conclude that $\tsr(B) \leq 2$.
\hfill$\qe$

\vskip 3mm

\begin{coro}
Let $A$ be a tsr boundedly divisible, unital C*-algebra with $\tsr(A) = 1$,
$G$ a finite group, and $\alpha$ an action of $G$ on $A$. Then
$A \rtimes_\alpha G$ is tsr boundedly divisible. Moreover
we have $\tsr(A \rtimes_\alpha G) \leq 2$.
\end{coro}

{\it Proof}.
Since the inclusion $1 \in A \subset A \rtimes_\alpha G$ is of
index-finite
type with depth 2 by Lemma~ 3.4, we can get the statement from Theorem 5.1.
\hfill$\qe$

\vskip 5mm

\begin{coro}
Let $A$ be a unital C*-algebra with $\tsr(A) =1$,
$G$ a finite group, and $\alpha$ an action of $G$ on $A$.
Then $(A \otimes UHF) \rtimes_{\alpha \otimes id} G$ is
tsr boundedly divisible.
\end{coro}

{\it Proof}.
Since $A \otimes UHF$ is tsr boundedly  divisible
and $\tsr(A \otimes UHF) = 1$ by using Theorem~ 2.3 (2),
it comes from Corollary~ 5.5.
\hfill$\qe$

\vskip 5mm

When $A$ in Corollary~ 5.6 is the trivial C*-algebra ${\mathbb C}$,
we can get an affirmative data for B. Blackadar's question.

\begin{coro}
Let $A$ be a UHF C*-algebra. Let $G$ be a finite group, and
$\alpha$ be an action of $G$ on $A$. Then
$A \rtimes_\alpha G$ is tsr boundedly divisible and
$$
\tsr(A \rtimes_\alpha G) \leq 2.
$$
\end{coro}

\vskip 3mm

\begin{rk}
The estimate in {\em Theorem~ 5.1} is best possible.
Indeed in \cite[Exam
ple 8.2.1]{bl3}
B. Blackadar constructed an symmetry action $\alpha$ on $CAR$
such that
$$
(C[0,1] \otimes CAR) \rtimes_{id \otimes \alpha} Z_2 \cong C[0,1] \otimes B,
$$
where $B$ is the Bunce-Deddens algebra of type $2^\infty$.
Then since $K_1(B)$ is non-trivial, we know that
$$
\tsr(C[0,1] \otimes B) = 2.
$$
$($See also \cite[Proposition~ 5.2]{nop}.$)$
\end{rk}

\vskip 5mm

Before ending this section we present an inclusion
$1 \in A \subset B$ with index $2$ such that $A$ is a tsr boundedly
divisible and $B$ can not be realized as some crossed product
algebra of $A$ by ${\mathbb Z}/2{\mathbb Z}$.

Let $A$ be a unital $C^*$-algebra, $\alpha$ an action of
${\mathbb Z}/2{\mathbb Z}$ on $A$, and
$B$ be the crossed product $A \rtimes_{\alpha} {\mathbb Z}/2{\mathbb Z}$.
Let $E$ the canonical conditional expectation and $u$
 an implemented unitary $u$ such that $\alpha(u) = u au^*$
for $a \in A$ as in the proof of Lemma~ 3.4.
Then $\{(1,1), (u^*, u)\}$ is a quasi-basis for $E$.
Note that $E(u) = 0$.
By \cite[Lemma~ 3.3.4]{wata} there is a *-isomorphism
$\varphi : C^*(B, e_A) \rightarrow qM_2(A)q$
such that
$$
q = \left(
\begin{array}{cc}
E(1\cdot 1) & E(u^*) \\
E(u) & E(uu^*)
\end{array}
\right)
=
\left(
\begin{array}{cc}
1 & 0 \\
0 & 1
\end{array}
\right)
$$
and
$$
\varphi(xe_Ay) = \left(\begin{array}{cc}
           E(x)E(y)&E(x)E(yu^*)\\
           E(ux)E(y)&E(ux)E(yu^*)
           \end{array}
           \right)
$$
for $x, y \in B$. Here $e_A$ is the Jones projection for the inclusion
$A \subset B$.
Therefore we can identify the basic construction with $M_2(A)$.

By this identification,
$$
A \cong \left\{
\left(
\begin{array}{cc}
a & 0 \\
0 & \alpha(a)
\end{array}
\right)
\ | \ a \in A
\right\}
\quad
B \cong \left\{
\left(
\begin{array}{cc}
a & b \\
\alpha(b) & \alpha(a)
\end{array}
\right)
\  | \  a, b \in A
\right\}
$$
and
$$
\varphi(e_A) =
\left(
\begin{array}{cc}
1 & 0 \\
0 & 0
\end{array}
\right), \quad
\varphi(1 - e_A) =
\left(
\begin{array}{cc}
0 & 0 \\
0 & 1
\end{array}
\right), \quad
\varphi(u
) =
\left(
\begin{array}{cc}
0 & 1 \\
1 & 0
\end{array}
\right).
$$
It is easy to see that $[\varphi(e_A)] = [1 - \varphi(e_A)]$  in $K_0(A)$.

{}From this observation we have

\vskip 3mm

\begin{ho}
Let $A$ be a unital C*-algebra and $\alpha$ an action of
${\mathbb Z}/2{\mathbb Z}$ on $A$. Let $B$ be the crossed product
$A \rtimes_\alpha{\mathbb Z}/2{\mathbb Z}$ and $e_A$ the Jones
projection of the inclusion $A \subset B$. If $\varphi : C^*(B,e_A)
\rightarrow M_n(A)$ is the canonical isomorphism, then
we have $[\varphi(e_A)] = [1 - \varphi(e_A)]$ in $K_0(A)$.
\end{ho}

\vskip 3mm

\begin{pro}
Let $A$ be a simple unital C*-algebra.
Suppose that $p$ is a projection in $A$ with $[p] \not= [\alpha(p)]$ in
$K_0(A)$.
Then the inclusion $pAp \subset pBp$ can not be represented as
$pAp \subset pAp \times_{\beta}{\mathbb Z}/2{\mathbb Z}$ for any
$\beta \in Aut(pAp)$.
\end{pro}

\vskip 3mm

{\it Proof}.
By the identification,
$
\varphi(p) = \left(
\begin{array}{cc}
p & 0 \\
0 & \alpha(p)
\end{array}
\right)
$.  The Jones projection for the conditional expectation $F_p = E|_{pBp}:
pBp
\to pAp$ is $e_Ap$ by Proposition~ 4.3, so
$\varphi(e_Ap) =  \left(
\begin{matrix}
p & 0 \\
0 & 0
\end{matrix}
\right)
$.
Then $\varphi(p -e_Ap) = \left(
\begin{matrix}
0 & 0 \\
0 & \alpha(p)
\end{matrix}
\right).
$
By the assumption, $[\varphi(e_Ap)] \not= [\varphi(p- e_Ap)]$ in $K_0(A)$
and hence
$[\varphi(e_Ap)] \not= [\varphi(p- e_Ap)]$ in $K_0(pAp)$.
So we have the conclusion by Lemma~ 5.9.
\hfill$\Box$

\vskip 3mm

\begin{rk}
When $A = A_1 \oplus A_2$ for simple unital C*-algebras $A_1$ and $A_2$,
we can also get the same conclusion in Proposition~ 5.10. Indeed,
$e_Ap$ becomes the Jones projection of the inclusion
$pAp \subset pBp$ as the same argument in Proposition 4.3.
\hfill$\Box$
\end{rk}

\vskip 3mm

The following example is due to T. Katsura and N. C. Phillips.

\begin{ex}
Let $\alpha$ be an automorphism on $CAR \otimes {\mathbb K}$ such that
$[\alpha(1 \otimes e_0)] \not= [1 \otimes e_0]$ in $K_0(CAR)$, where
$e_0$ is a minimal projection. Such an automorphism can be constructed
by modifying the shift operator in \cite{ki2}.
Set $D$ as the unitaization of $CAR \otimes {\mathbb K}$.
Then $\alpha$ can be an automorphism on $D$. We call it $\alpha$ again.

Define a symmetry $\gamma$ on $D \oplus D$ by
$\gamma((a,b)) = (\alpha^{-1}(b),\alpha(a))$ for $(a,b)$ in $D \oplus D$.
Consider the inclusion
$$
D \oplus D \subset (D \oplus D) \rtimes_\gamma{\mathbb Z}/2{\mathbb Z}.
$$

Since
$$
\begin{array}{ll}
[\gamma((1 \otimes e_0, 1 \otimes e_0))]
&= [(\alpha^{-1}(1 \otimes e_0), \alpha(1 \otimes e_0)]\\
&\not= [(1 \otimes e_0, 1 \otimes e_0)]
\end{array}
$$
in $K_0(D \oplus D)$ by the construction. We know that the inclusion
$$
(1 \otimes e_0, 1 \otimes e_0)(D \oplus D) (1 \otimes e_0, 1 \otimes e_0)
\subset (1 \otimes e_0, 1 \otimes e_0)
((D \oplus D) \rtimes_\gamma{\mathbb Z}/2{\mathbb Z})
(1 \otimes e_0, 1 \otimes e_0)
$$
can not be represented as
$$
(1 \otimes e_0, 1 \otimes e_0)(D
 \oplus D)(1 \otimes e_0, 1 \otimes e_0)
\subset
(1 \otimes e_0, 1 \otimes e_0)(D \oplus D)(1 \otimes e_0, 1 \otimes e_0)
\rtimes_\beta {\mathbb Z}/2{\mathbb Z}
$$
for any $\beta \in
Aut((1 \otimes e_0, 1 \otimes e_0)(D \oplus D)(1 \otimes e_0, 1 \otimes
e_0))$
by Proposition~ 5.10  and Remark~ 5.11.
Note that
$$
(1 \otimes e_0, 1 \otimes e_0)(D \oplus D)(1 \otimes e_0, 1 \otimes e_0)
\cong CAR \oplus CAR,
$$ that is, the algebra is tsr boundedly divisible.
\hfill$\Box$
\end{ex}


\vskip 5mm


$$
\begin{array}{ll}
\hbox{HIroyuki Osaka} & \hbox{Tamotsu Teruya}\\
\hbox{Department of Mathematics} &\hbox{Department of Mathematical Sciences}\\
\hbox{University of Oregon} &\hbox{Ryukyu University}\\
\hbox{Eugene, Oregon 97403-1222} &\hbox{Nishihara-cho, Okinawa 903-0213}\\
\hbox{U.S.A}. &\hbox{Japan}\\
&\\
\hbox{current address} & \hbox{hebo2@nirai.ne.jp}\\
&\\
\hbox{Department of Mathematical Sciences} &\\
\hbox{Ritsumeikan University} &\\
\hbox{Kusatsu, Shiga 525-8577} &\\
\hbox{Japan} & \\
\hbox{osaka@se.ritsumei.ac.jp} &
\end{array}
$$

\end{document}